\documentclass[12pt]{article}

\usepackage{algorithm, algpseudocode}
\usepackage{graphicx}
\usepackage{float}
\usepackage{tikz}
\usepackage{capt-of}
\usepackage{color}
\usepackage{fullpage}
\usepackage{amssymb}
\usepackage{amsmath}
\usepackage{amsthm}
\usepackage{amsfonts}
\usepackage{enumitem}
\usepackage{amscd}

\usepackage[colorlinks=true,
linkcolor=brown,
filecolor=brown,
citecolor=brown]{hyperref}

\newcommand{\seqnum}[1]{\href{https://oeis.org/#1}{\rm \underline{#1}}}

\DeclareMathOperator{\Sum}{Sum}
\DeclareMathOperator{\SumSq}{SumSq}

\begin{document}

\theoremstyle{definition}
\newtheorem{theorem}{Theorem}
\newtheorem{corollary}[theorem]{Corollary}
\newtheorem{lemma}[theorem]{Lemma}
\newtheorem{proposition}[theorem]{Proposition}
\newtheorem{definition}[theorem]{Definition}
\newtheorem{example}[theorem]{Example}

\begin{center}
\vskip 1cm{\LARGE\bf 
Counting Interval Sizes in the Poset of\\ Monotone Boolean Functions
}
\vskip 1cm
\large
Bart\l{}omiej Pawelski\\
Institute of Informatics\\
University of Gda\'nsk\\
\href{mailto:bartlomiej.pawelski@ug.edu.pl}{\tt bartlomiej.pawelski@ug.edu.pl} \\
\end{center}

\vskip .2in

\begin{abstract}
We focus on the computational aspects of counting interval sizes in the poset $D_n$, which represents all monotone Boolean functions of $n$ variables. We present a resource-aware algorithm enabling the calculation of interval sizes in $D_7$.
\end{abstract}
    
\section{Introduction}\label{sec:intro}

Consider the set $B = \{0,1\}$. For any integer $n$, the set $B^n$ stands for the set of $n$-element sequences of $B$. A Boolean function is a mapping $f: B^n \to B$. The $B^n$ contains $2^n$ members, and there are $2^{2^n}$ Boolean functions of $n$ variables.

There is a partial order in $B$: $0 \leq 1$, which extends naturally to $B^n$. Specifically, for any pair of elements $x= (x_1,\dots,x_n)$ and $y=(y_1,\dots,y_n)$ from $B^n$, we say $x\le y$ if for each $i$, the condition $x_i \leq y_i$ is satisfied.

A Boolean function is said to be \textit{monotone} if, when $x \leq y$, then $f(x) \leq f(y)$. We define $D_n$ as the set of all monotone Boolean functions of $n$ variables. Let $d_n$ represent the cardinality of $D_n$, which is also known as the $n$--th Dedekind number. Dedekind numbers are described by the \textit{On-Line Encyclopedia of Integer Sequences} (OEIS) sequence \seqnum{A000372}. Dedekind numbers hold significance in various combinatorial contexts, as they correspond to the number of simple games with $n$ players in minimal winning form, the number of antichains of subsets of an $n$ set, the number of Sperner families and the cardinality of a free distributive lattice on $n$ generators \cite{markowsky, oeis}.

A partial order in $D_n$ is defined as follows: for two functions $f, g \in D_n$, $f \leq g$ if condition $f(x) \leq g(x)$ is satisfied for every $x \in B^n$. We denote an \textit{interval} $[f, g]$ as the set of functions $h \in D_n$ that satisfy $f \leq h \leq g$. Let $\#[f, g]$ denote the cardinality of $[f, g]$. Let $\#[f, \top]$ denote the number of functions $g$ that satisfy $f \leq g$.

Two monotone Boolean functions are said to be \textit{equivalent} if one function can be obtained from the other by permuting the input variables. For $h \in D_n$, let $\gamma(h)$ denote the number of functions in $D_n$ equivalent to $h$. Let $R_n$ denote the set of equivalence classes of $D_n$, with $r_n$ representing its cardinality. This cardinality is also referred to as the number of inequivalent monotone Boolean functions of $n$ variables. The values of $r_n$ are listed in the OEIS sequence \seqnum{A003182}.

We consider the problem of counting the interval size $\#[f, g]$ for $f, g \in D_n$. The most recent advancement in this topic is the computation of the sizes of intervals $\#[f, \top]$ for all $f \in D_7$, as achieved by Van Hirtum \cite{hirtum}. Those sizes were used to calculate $d_9$ by Van Hirtum, De Causmaecker, Goemaere, Kenter, Riebler, Lass, and Plessl \cite{hirtum2}. The value of $d_9$ was independently calculated by Jäkel \cite{jakel}. In this paper, we introduce a novel algorithm that accelerates this computation process of interval sizes even further. With our implementation of this algorithm, we can calculate $\#[f, \top]$ for all $f \in D_7$ in approximately 2.5 hours.

\begin{table*}
\begin{center} \renewcommand{\arraystretch}{1}
\resizebox{\linewidth}{!}{%
\begin{tabular}{l|l|l} $n$ & $d_n$ & $r_n$\\ \hline 
0&  2  & 2    \\  
1&  3  & 3    \\
2&  6  & 5     \\  
3&  20 & 10     \\  
4&  168 & 30      \\
5&  7581 & 210      \\  
6&  7828354 & 16353    \\  
7&  2414682040998 & 490013048 \\
8&  56130437228687557907788 & 1392195548889993358 \\
9&  286386577668298411128469151667598498812366 & 789204635842035040527740846300252680 \\ \hline \end{tabular}} \caption{Known values
of $d_n$ and $r_n$.} \label{tab:valuesd9}
\end{center}
\end{table*}

\section{Methodology}

\subsection{Posets}

A binary relation that is reflexive, antisymmetric and transitive, when defined on a set $P$, forms a \textit{partially ordered set}, or simply a \textit{poset}. In the Introduction, we defined partial orders on $B$, $B^n$ and $D_n$, making them posets.

Two posets $P = (X, \leq)$ and $Q = (Y, \leq)$ are said to be \textit{isomorphic} if there exists a bijection $f: X \to Y$ so that for any elements $x_1, x_2$ in $P$, $x_1 \leq x_2$ holds true iff $f(x_1) \leq f(x_2)$. We use the equals symbol $=$ to represent an isomorphism between two posets.

For a given poset $P = (X, \le)$, the \textit{incidence matrix} $M_P$ is a binary matrix where the rows and columns correspond to elements of $X$. An element $M_P(x,y)$ in the matrix is 1 if $x \le y$ and 0 otherwise.

Given two posets, $(X, \le)$ and $(Y, \le)$, their Cartesian product is represented as $X \times Y$. In this product, $(a,b) \le (c,d)$ holds true iff $a \le c$ and $b \le d$.
A function $f: X \to Y$ is monotone when, for any elements $x,y \in X$ with $x \leq y$, we have $f(x) \leq f(y)$. We represent the set of all monotone functions from $X$ to $Y$ as $Y^X$. When considering two functions $f,g \in Y^X$, we define $f \leq g$ to be true if $f(x) \leq g(x)$ for every $x \in X$. Following this notion, we have $D_n = B^{B^n}$.

\begin{lemma}\label{L2} \cite{szepietowski}
\leavevmode
\begin{itemize}
\item[(a)] \quad $B^{k+m}=B^{k}\times B^{m}$.
\item[(b)] \quad $D_{k+m}=(D_k)^{B^m}$.
\end{itemize}
\end{lemma}

This lemma, presented in different formulations, is commonly used in literature for computing Dedekind numbers. Wiedemann \cite{Wied1991} used the isomorphism $D_8 = (D_6)^{B^2}$ to compute $d_8$, and J\"akel \cite{jakel} used the isomorphism $D_9 = (D_5)^{B^4}$ to compute $d_9$.

\subsection{Monotone Boolean functions}

We represent any Boolean function of $n$ variables using a binary sequence of length $2^n$, termed its truth table.

\begin{example}
   Consider the following truth table:

   \begin{table}[H]
\centering
{\renewcommand{\arraystretch}{1.2}
\begin{tabular}{c|c|c|c|c|c|c|c}
$000$ & $001$ & $010$ & $011$ & $100$ & $101$ & $110$ & $111$  \\ \hline
0 & 1 & 1 & 1 & 0 & 1 & 1 & 1 \\
\end{tabular}}
\end{table}

The first row contains lexicographically ordered elements of $B^3$, and the second row contains a function value $f(x)$ for a given $x \in B^3$. We represent this Boolean function by the binary word 01110111.

\end{example} 

Note that functions with six or fewer variables can be conveniently represented as a 64-bit integer, which corresponds to a word length of modern CPUs, making it a notable efficient representation for computational purposes.

For $f, g$ $\in D_n$, the union of these functions is represented by $f \cup g$, and the intersection by $f \cap g$. In terms of their binary representations, the union can be represented using the bitwise OR operation and the intersection using the bitwise AND operation.

According to Lemma \ref{L2}, each function $f \in D_{n+1}$ can be represented as a concatenation of two functions $(f_0, f_1)$, where $f_0, f_1 \in D_n$, and $f_0 \leq f_1$. This lemma also implies that any function $f \in D_{n+m}$ can be decomposed into $2^m$ functions from $D_n$.

\begin{example}
    The two functions in $D_0$ are: 0 and 1. The three functions in $D_1$ are: 00, 01, and 11. The six functions in $D_2$ are: 0000, 0001, 0011, 0101, 0111, and 1111.
\end{example}

\begin{example}
Consider the monotone Boolean function $f$ represented by the binary sequence 00110111, which belongs to $D_3 = B^{B^3}$. The function $f$ can be decomposed into a tuple of two functions $(0011, 0111)$, where each component is an element of $D_2$ and 0011 $\leq$ 0111.

\end{example}

\section{Algorithms}

In this chapter, we explore two approaches designed to calculate the interval sizes within posets, specifically focusing on $D_n$. First algorithm is well--known and can be found in the literature (for example see: \cite{aigner, Fid2001}). The second approach is specifically designed to address the problem of counting interval sizes in $D_7$.
\subsection{Finding interval sizes in $D_n$}

For matrix $M$, we define:

\[
\Sum(M) = \sum_{i,j} M(i, j),
\]
\[
\SumSq(M) = \sum_{i,j} (M(i, j))^2.
\]

Let us recall that we use $M_{D_n}$ to denote the incidence matrix of $D_n$. By $(M_{D_n})^2$ we denote the product $M_{D_n} \times M_{D_n}$. The core of the following approach comes from the known results in the incidence algebra.

\begin{proposition} \label{lemInc} \cite{aigner}[Proposition 4.7]
    For any $x, y \in D_n$, we have $\#[x, y] = (M_{D_n})^2(x,y)$.
\end{proposition}

\begin{example}
Consider $M_{D_2}$. Let us recall $D_2 = \{0000, 0001, 0011, 0101, 0111, 1111\}$:
$$
M_{D_2}=
\begin{bmatrix}
1& 1& 1& 1& 1& 1 \\
0& 1& 1& 1& 1& 1 \\
0& 0& 1& 0& 1& 1 \\
0& 0& 0& 1& 1& 1 \\
0& 0& 0& 0& 1& 1 \\
0& 0& 0& 0& 0& 1
\end{bmatrix}.
$$

Consider the square of $M_{D_2}$: 
$$
(M_{D_2})^2=
\begin{bmatrix}
1& 2& 3& 3& 5& 6 \\
0& 1& 2& 2& 4& 5 \\
0& 0& 1& 0& 2& 3 \\
0& 0& 0& 1& 2& 3 \\
0& 0& 0& 0& 1& 2 \\
0& 0& 0& 0& 0& 1
\end{bmatrix}.
$$

From Proposition \ref{lemInc}, the matrix $(M_{D_2})^2$ contains all interval sizes in $D_n$. Note that one can utilize a matrix $(M_{D_n})^2$ to calculate $d_{n+2}$ \cite{Fid2001}[Algorithm 2.]:

$$d_{n+2} = \SumSq((M_{D_n})^2).$$

\end{example}

Unfortunately, in practice, this method cannot be applied for $D_n$ with $n$ greater than 5 due to the size of such a data structure. For $n=5$, the incidence matrix in the form of $7581 \times 7581$ bits should occupy approximately 58 megabytes (using bitsets) and would require double that amount of storage if the columns are stored separately as well. However, for $n=6$, the incidence matrix would consist of $7828354 \times 7828354$ elements, which would require around 62 terabytes of storage. This is far beyond what is achievable with the storage capabilities of standard personal computers.

\subsection{Finding interval sizes in $D_{n+2}$}

In the second algorithm, we focus on avoiding storing a whole incidence matrix of $D_n$. Instead of this, we use the isomorphism:

$$D_{n+2} = D_{n} ^{B^{2}}.$$

Consider the function $x \in D_{n+2}$. It can be decomposed into $x = (x_0, x_1, x_2, x_3)$ with $x_i \in D_n$. The tuple $y = (y_0, y_1, y_2, y_3)$ with $y_i \in D_n$ belongs to $[x, \top]$ if and only if:

\begin{enumerate}[label=]
    \item $y_0 \in [x_0, \top],$
    \item $y_3 \in [x_3, \top],$
    \item $y_1 \in [y_0 \cup x_1, y_3],$
    \item $y_2 \in [y_0 \cup x_2, y_3].$
\end{enumerate}

Thus, we have:

$$\#[x, \top] = \sum_{y_0 \in [x_0, \top]} \sum_{y_3 \in [x_3, \top]} (M_{D_n})^2((y_0 \cup x_1), y_3) \cdot (M_{D_n})^2((y_0 \cup x_2), y_3).$$

Using this approach, we can construct an algorithm that finds the size of the interval $\#[x, \top]$ for any $x \in D_{n+2}$. Note that the algorithm uses the precomputed matrix $(M_{D_n})^2$.

\begin{algorithm}[H]
\caption{Calculation of $\#[x, \top]$ for $x \in D_{n+2}$}
    \vspace{1mm}
    \hspace*{\algorithmicindent} \textbf{Input:} $x \in D_{n+2}$ \\
    \hspace*{\algorithmicindent} \textbf{Output:} $s = \#[x, \top]$
\begin{algorithmic}[1]
\State Initialize $s = 0$,
\State Decompose $x \in D_{n+2}$ into $x_0, x_1, x_2, x_3 \in D_{n}$,
\ForAll {$y_0 \in [x_0, \top] $}
\ForAll {$y_3 \in [x_3, \top] $}
\State $s = s + (M_{D_n})^2((y_0 \cup x_1), y_3) \cdot (M_{D_n})^2((y_0 \cup x_2), y_3)$
\EndFor
\EndFor
\end{algorithmic}
\end{algorithm}

Similarily, we can construct an algorithm finding the size of interval $\#[x, y]$ for any $x, y \in D_{n+2}$.

\begin{algorithm}[H]
\caption{Calculation of $\#[x, y]$ for $x, y \in D_{n+2}$ }
    \vspace{1mm}
    \hspace*{\algorithmicindent} \textbf{Input:} $x, y \in D_{n+2}$ \\
    \hspace*{\algorithmicindent} \textbf{Output:} $s = \#[x, y]$
\begin{algorithmic}[1]
\State Initialize $s = 0$,
\State Decompose $x \in D_{n+2}$ into $x_0, x_1, x_2, x_3 \in D_{n}$,
\State Decompose $y \in D_{n+2}$ into $y_0, y_1, y_2, y_3 \in D_{n}$,
\ForAll {$f_0 \in [x_0, y_0] $}
\ForAll {$f_3 \in [x_3, y_3] $}
\State $s = s + (M_{D_n})^2((f_0 \cup x_1), (f_3 \cap y_1)) \cdot (M_{D_n})^2((f_0 \cup x_2), (f_3 \cap y_2))$ 
\EndFor 
\EndFor
\end{algorithmic}
\end{algorithm}

\section{Implementation}

According to the Introduction, our objective is to calculate the interval size $\#[x, \top]$ for all $x \in D_7$. This feat was previously accomplished only by Van Hirtum in 2021 \cite{hirtum}.

To accelerate computations, we utilize the property that for every $x \in D_n$ and every permutation $\pi \in S_n$ (see \cite{Wied1991}):
$$
    \# [x, \top] = \# [\pi(x), \top].
$$

Consequently, instead of calculating $\# [x, \top]$ for all 2,414,682,040,998 functions from $D_7$, we can limit our computations to just 490,013,148 functions from $R_7$, achieving a speed--up by a factor of nearly 5,000.

We calculated $\# [x, \top]$ for all $x \in R_7$ using our implementation of Algorithm 1 in approximately 2.5 hours on a 32-thread Xeon machine. This demonstrates a significant improvement in efficiency compared to previous methods; for instance, a similar computation conducted by Van Hirtum required 6.5 hours on a 36-core computer.

In order to confirm our results, we performed the following test:

$$
    d_8 = \sum_{x\in R_7} \#[x, \top] \cdot \gamma(x) = 56130437228687557907788,
$$

which is equal to the actual value of $d_8$ (see Table \ref{tab:valuesd9}).


\newpage

\begin{thebibliography}{99}

    \bibitem{aigner} M. Aigner.
    \newblock {\it Combinatorial Theory.} 
    \newblock Springer, 1979.

    \bibitem{Fid2001}
    R. Fidytek, A. W. Mostowski, R. Somla, A. Szepietowski.
    \newblock Algorithms counting monotone Boolean functions
    \newblock in {\em Information Processing Letters}, 2001, issue 79, pages 203--209

    \bibitem{hirtum}
    L. Van Hirtum.
    \newblock A path to compute the 9th Dedekind
    Number using FPGA Supercomputing (2021). Master's thesis.
    \newblock {\url{https://hirtum.com/thesis.pdf}}, accessed 01.11.2023

    \bibitem{hirtum2}
    L. Van Hirtum, P. De Causmaecker, J. Goemaere, T. Kenter, H. Riebler, M. Lass, and C. Plessl, A computation of D(9) using FPGA supercomputing, preprint, 2023. \\ Available at \url{https://arxiv.org/abs/2304.03039}.

    \bibitem{jakel}
    C. J\"akel, A computation of the ninth Dedekind number, preprint, 2023. \\ Available at \url{https://arxiv.org/abs/2304.00895}.

    \bibitem{markowsky}
    G. Markowsky, Free completely distributive lattices, \emph{Proc. Amer. Math. Soc.} {\bf 74} (1979).
    
    \bibitem{pawelski}
    B. Pawelski, On the number of inequivalent monotone Boolean functions of 8 variables,
    \emph{J. Integer Sequences} {\bf 25} (2022),
    \href{https://cs.uwaterloo.ca/journals/JIS/VOL25/Pawelski/pawelski7.html}{Article 25.7.7}.

    \bibitem{szepietowski}
    A. Szepietowski, Fixes of permutations acting on monotone Boolean functions,
    \emph{J. Integer Sequences} {\bf 25} (2022),
    \href{https://cs.uwaterloo.ca/journals/JIS/VOL25/Szepietowski/szep7.html}{Article 25.9.6}.

    \bibitem{oeis}
    OEIS Foundation Inc. (2023).
    \newblock Dedekind numbers, Entry A000372 in The On-Line Encyclopedia of Integer Sequences, {\url{https://oeis.org/A000372}}
    

    \bibitem{Wied1991}
    D. Wiedemann.
    \newblock A computation of the eighth Dedekind number
    \newblock in {\em Order 8}, 1991, pages 5--6
    
    \bibitem{yamamoto}
    K. Yamamoto.
    \newblock Note on the Order of Free Distributive Lattices
    \newblock in {\em The Science Reports of the Kanazawa University}, 1953, Volume II, No. 1, pages 5--6
    \end{thebibliography}
\end{document}